\documentclass[12pt]{amsart}
\usepackage{amsmath,amsthm}
\usepackage{mathrsfs}
\input amssym.def
\input amssym
\usepackage{enumerate}
\usepackage{hyperref}
\hypersetup{
    pdftitle={Super},    % title
    pdfauthor={Brendan Miller}{Ilya Krishtal},     % author
    pdfsubject={Kadec for DF},   % subject of the document
    pdfcreator={Ilya Krishtal},   % creator of the document
    pdfproducer={Ilya Krishtal}, % producer of the document
    pdfkeywords={Banach modules}{Beurling spectrum}{Kadec theorem}, % list of keywords
     colorlinks,%
    citecolor=blue,%
    filecolor=black,%
    linkcolor=magenta,%
    urlcolor=black
}

\usepackage{algorithm}
\usepackage{algorithmicx}
\usepackage[noend]{algpseudocode}
\usepackage{tikz}
\allowdisplaybreaks

\chardef\bslash=`\\ % p. 424, TeXbook
\makeatletter
\def\verbatim{\interlinepenalty\@M \@verbatim
   \leftskip\@totalleftmargin\advance\leftskip2pc
   \frenchspacing\@vobeyspaces \@xverbatim}
\makeatother
\newtheorem{thm}{Theorem}[section]
\newtheorem{cor}[thm]{Corollary}
\newtheorem{lem}[thm]{Lemma}
\newtheorem{prop}[thm]{Proposition}

\theoremstyle{definition}
\newtheorem{defn}{Definition}[section]

\theoremstyle{remark}
\newtheorem{rem}{Remark}[section]
\newtheorem{exmp}{Example}[section]

\numberwithin{equation}{section}

\newcommand{\begeq}{\begin {equation}}
\newcommand{\eq}{\end{equation}}
\newcommand{\bs}{\begin {split}}
\newcommand{\es}{\end{split}}
\newcommand{\bp}{\begin {prop}}
\newcommand{\ep}{\end {prop}}
\newcommand{\bt}{\begin {thm}}
\newcommand{\et}{\end {thm}}
\newcommand{\bc}{\begin {cor}}
\newcommand{\ec}{\end {cor}}
\newcommand{\bl}{\begin {lem}}
\newcommand{\el}{\end {lem}}
\newcommand{\bpf}{\begin {proof}}
\newcommand{\epf}{\end {proof}}
\newcommand{\bi}{\begin {itemize}}
\newcommand{\ei}{\end {itemize}}
\newcommand{\ben}{\begin {enumerate}}
\newcommand{\een}{\end {enumerate}}
\newcommand{\brem}{\begin {rem}}
\newcommand{\erem}{\end {rem}}

\newcommand{\bd}{\begin {defn}}
\newcommand{\ed}{\end {defn}}

\newcommand{\bex}{\begin {exmp}}
\newcommand{\eex}{\end {exmp}}

\newcommand{\norm}[1]{\left\|{#1} \right\|}

\newcommand{\la}{\langle}
\newcommand{\ra}{\rangle}

\newcommand{\A}{\mathcal A}

\newcommand{\F}{\mathcal{F}}

\newcommand{\TTT}{{\T\kern-.44em \T}}
\newcommand{\tTTT}{\widetilde{\T\kern-.44em \T}}
\newcommand{\chT}{\check\T}

\newcommand{\ZZ}{{\mathbb Z}}

\newcommand{\RR}{{\mathbb R}}
\newcommand{\CC}{{\mathbb C}}

\newcommand{\X}{\mathcal{X}}

\newcommand{\HH}{{\mathcal H}}

\renewcommand{\d}{\delta}

\renewcommand{\l}{\lambda}
\renewcommand{\L}{\Lambda}

\newcommand{\lhat}{\mathcal FL^1}%{\widehat{L_{1}}}

\newcommand{\lloc}{\mathcal FL^{1}_{loc}}

\newcommand{\g}{\gamma}

\newcommand{\T}{\mathcal{T}}

\DeclareMathOperator{\supp}{supp}

\definecolor{bmcolor}{rgb}{0.9, 0.3, 0}

\definecolor{ikcolor}{rgb}{1, 0, 1.}

\begin{document}

\title {Kadec-type theorems for sampled group orbits}

%%%%%%%%%%%%%%%%%%%%%%%%%%%%%%%%%%%%%%%%%%%%%%%%%%%%%%%%%%%%%%%%%%%%%%%%
\author{ Ilya Krishtal and Brendan Miller }
\address{Department of Mathematical Sciences, Northern Illinois University, DeKalb, IL 60115 \\
email: ikrishtal@niu.edu, bmiller14@niu.edu}
 
%\thanks{ The first author is supported in part by RFBR grant 13-01-00378.
%The second author is supported in part by NSF grant DMS-0908239.}

%%%%%%%%%%%%%%%%%%%%%%%%%%%%%%%%%%%%%%%%%%%%%%%%%%%%%%%%%%%%%%%%%%%%%%%%

\date{\today }

%\subjclass[2010]{46H25}

\keywords{Perturbation of bases and frames, Banach modules, Beurling spectrum}

\dedicatory{to Anatoly G.~Baskakov on the occasion of his $80^{ th}$ birthday.}

%%%%%%%%%%%%%%%%%%%%%%%%%%%%%%%%%%%%%%%%%%%%%%%%%%%%%%%%%%%%%%%%%%%%%%%%
%%%%%%%%%%%%%%%%%%%%%%%%%%%%%%%%%%%%%%%
%%%%%%%%%%%%%%%%%%%%%%%%%%%%%%%%%
\begin{abstract}
We extend the classical  Kadec $\frac14$ theorem for systems of exponential functions on an interval to frames and atomic decompositions formed by sampling an orbit of a vector under an isometric group representation. 
\end{abstract}
%%%%%%%%%%%%%%%%%%%%%%%%%%%%%%%%%%%%%%%
%%%%%%%%%%%%%%%%%%%%%%%%%%%%%%%%%
\maketitle

\section{Introduction and main results}

The celebrated Kadec $\frac14$ theorem \cite{K64} states that a system of exponential functions $\{\mathbf{e}_{\mu_n}: n\in\ZZ\}$, where $\mathbf{e}_{\mu_n}(t) = e^{i\mu_n t}$, $t\in (-\pi, \pi)$, is a Riesz basis (see Definition \ref{defRB}) for the Hilbert space $\HH = L^2(-\pi, \pi)$ provided that $\delta = \sup\{|\mu_n-n|: n\in\ZZ\} <\frac14$. The purpose of this note is to show that this type of perturbation result holds not only for exponential bases but in a much more general setting where the basis is a sampled orbit of a vector under an isometric group representation.
Indeed, we have 
    \[
    \mathbf{e}_{\mu_n} = M(\mu_n)\mathbf{e}_0,  
    \]
where $M: \RR\to B(\HH)$, $(M(t)x)(s) = e^{its}x(s)$, $x\in\HH$, $s,t\in\RR$, is the modulation representation of the group $\RR$ by operators in $B(\HH)$ -- the space of all bounded linear operators on $\HH$. The main result of this paper is the following theorem which establishes that one can replace $M$ with any isometric representation $\T:\RR\to B(\HH)$ and $\mathbf{e}_0$ with any vector $x\in \HH$ such that its sampled orbit under $\T$ is a Riesz basis, on condition that the Beurling spectrum $\Lambda(\HH,\T)$, see Definition \ref{defBS}, is compact. In fact, following R.~Balan \cite{Bal97}, we establish that the result holds not just for Riesz bases but for frames as well, see Definition \ref{defFr}. 

\bt\label{thmRB}
Let $\T:\RR\to B(\HH)$ be an isometric representation such that $\Lambda(\HH,\T) \subseteq [-\g,\g]$ for some $\g > 0$. Assume that a vector $x\in\HH$ and a set $\Gamma = \{\g_n: n\in\ZZ\} \subset\RR$ are such that the system of vectors $\{\T(\g_n)x\}$ forms a frame for $\HH$ with bounds $B\ge A>0$. Let  $\widetilde\Gamma = \{\widetilde\g_n: n\in\ZZ\} \subset\RR$ be such that
\begin{equation}\label{deq}
    \delta := \sup\{|\widetilde\g_n-\g_n|: n\in\ZZ\} < \frac{\pi}{4\g} - \frac1\g\arcsin{\left(\frac{1}{\sqrt2}\left(1-\sqrt{\frac AB}\right)\right)}.
\end{equation}
Then the system of vectors $\{\T(\widetilde\g_n)x\}$ also forms a frame for $\HH$ with bounds
\[
A\left(1-\sqrt{\frac AB}(1-\cos\delta\g+\sin\delta\g)\right)^2 \mbox{ and } B\left(2-\cos\delta\g+\sin\delta\g\right)^2.
\]
If additionally $\{\T(\g_n)x\}$ is a Riesz basis for $\HH$ then so is $\{\T(\widetilde\g_n)x\}$.
\et

%\bt\label{thmFr}
%Let $\T:\RR\to B(\HH)$ be an isometric representation such that $\Lambda(\HH,\T) \subseteq [-\g,\g]$ for some $\g > 0$. Assume that a vector $x\in\HH$ and a separated set $\Gamma = \{\g_n: n\in\ZZ\} \subset\RR$ are such that the system of vectors $\{\T(\g_n)x\}$ forms a frame for $\HH$ with bounds $B\ge A>0$. Let  $\widetilde\Gamma = \{\widetilde\g_n: n\in\ZZ\} \subset\RR$ be such that inequality \eqref{deq} holds.
%\[
%\delta := \sup\{|\widetilde\g_n-\g_n|: n\in\ZZ\} < \frac{\pi}{4\g} - \frac1\g\arcsin{\left(\frac{1}{\sqrt2}\left(1-\sqrt{\frac AB}\right)\right)}.
%\]
%Then the system of vectors $\{\T(\widetilde\g_n)x\}$ forms a frame  for $\HH$ 
%\et

Observe that in the Riesz basis setting in the above theorem the result implies that the set $\Gamma$ has to be separated by $2\delta$ since the points in $\widetilde\Gamma$ may not be allowed to collide. Thus, the following corollary is immediate after applying standard trigonometric identities.

\bc
Assume that a vector $x\in\HH$ and a set $\Gamma = \{\g_n: n\in\ZZ\} \subset\RR$ are such that the system of vectors $\{\T(\g_n)x\}$ forms a Riesz basis for $\HH$ with bounds $B\ge A>0$. Then $\kappa := \inf\{|\g_m -\g_n|: m,n\in\ZZ\}$ satisfies
\[
2\left(1-\cos\frac{\g\kappa}{2}\right)\left(1+\sin\frac{\g\kappa}{2}\right)\ge\frac AB \quad \mbox{or}\quad \kappa \ge \frac{\pi}{2\g}.
\]
\ec

The key tool for proving Theorem \ref{thmRB} is the $\mathcal FL^1_{loc}$ operator functional calculus for generators of Banach $L^1(\RR)$-modules developed in \cite{BK14, BKU20}. We shall present necessary definitions and notation for it in the following section. The proof of Theorem \ref{thmRB} is then given in Section \ref{mainproof}. Since our tools remain valid in a Banach space setting, we also provide an extension of Theorerm \ref{thmRB} to atomic decompositions in Section \ref{adem1}.

Using the $\mathcal FL^1_{loc}$ functional calculus, one can obtain estimates for the operator norm $\|\T(t) - I\|$ for small values of $t> 0$, which is, in a way, a key ingredient for proving Kadec-type results. In \cite{B79}, see also \cite[Theorem 3.7]{BK05}, A.~Baskakov obtained the following estimate in a Banach module $(\X, \T)$ with $\Lambda(\X,\T) \subseteq [0,\gamma]$:
\begin{equation}\label{berneq}
    \|\T(t) - I\| \le 2\sqrt{2} \sup_{s\in [0,\gamma]} |e^{its} - 1| = 
4\sqrt2\sin\frac{\g t}{2}, \ t\in [0, \frac\pi\g].
\end{equation}
Clearly, it is only meaningful when $t \le \frac2\g\arcsin\frac1{2\sqrt2}$ as one always has $\|\T(t) - I\| \le 2$ for an isometric representation.

We conclude our introduction with the following significant improvement of the above inequality \eqref{berneq} which is in some sense similar to the improvement Kadec had to obtain on the way to proving his result.

\bt\label{Bineq}
Let $\T:\RR\to B(\X)$ be an isometric representation of the group $\RR$ by operators on a (complex) Banach space $\X$ such that $\Lambda(\X,\T) \subseteq [0,\g]$. Then for $t\in[0,\frac{\pi}{2\g}]$ one has 
\[
\|\T(t) - I\| \le 1 - \cos{\g t} + \sin{\g t} = 2\sqrt2\sin\frac{\g t}{2}\sin\left(\frac{\g t}{2}+\frac\pi4\right).
\]
\et

\section{$\mathcal F{L^1_{loc}}$ functional calculus}\label{calc}

In our presentation of the operator functional calculus for generators of Banach $L^1(\RR)$-modules we follow \cite{BK14, BKU20}.  
Let us recall some notation.

We denote by $\X$  a complex Banach space and by $B(\X)$ -- the Banach algebra of all bounded linear operators in $\X$. We  assume that $\X$ is endowed with a non-degenerate Banach module structure over the group algebra $L^1(\RR)$, which is associated with a strongly continuous representation $\T:\RR\to B(\X)$  of the locally compact Abelian group $\RR$ by operators in $B(\X)$. %The multiplication in $L^1(\RR)$ is the convolution 
%\[
%(f*g)(t) = \int_\RR f(s)g(t-s)ds, \ f,g\in L^1(\RR),\ t\in\RR.
%\]
\bd\label{dmod}
A complex Banach space $\X$  is a \emph{Banach module} over $L^1(\RR)$ if there is a bilinear map
$(f,x)\mapsto fx: L^1(\RR)\times\mathcal X\to\mathcal X$ which has the following properties:
\begin{enumerate}
\item $(f*g)x = f(gx)$,  $f,g\in L^1(\RR)$, $x\in\mathcal X$;
\item $\norm{fx} \le \|f\|_1\|x\|$,  $f\in L^1(\RR)$, $x\in\mathcal X$.
\end{enumerate}
\ed
As usual (see, e.g., \cite{BK05}), by
non-degeneracy of the module we mean that $fx = 0$ for all $f\in L^1(\RR)$ implies that $x = 0$. The condition that the module structure is associated with a strongly continuous representation $\T$ means that it is defined via the Bochner integral
\begin{equation}\label{assrep}
  fx = \int_\RR f(t)\T(-t)xdt, f\in L^1(\RR), x\in\X.  
\end{equation}
 With a slight abuse of notation \cite{BK05}, given $f\in L^1(\RR)$, we shall denote by $\T(f)$ the operator in $B(\X)$ defined by $\T(f)x = fx$, $x\in\X$. Observe that we have $\|\T(f)\|\le \|f\|_1$, $f\in L^1(\RR)$, by Property (2) in Definition \ref{dmod}. For a Banach module $\X$, we will also use the notation $(\X,\T)$ if we want to emphasize that the module structure is associated with the representation $\T$.

We use the Fourier transform in the form
\[
(\F(f))(\xi) = \widehat f(\xi) = \int_\RR f(t)e^{-it\xi}dt, \ f\in L^1(\RR).
\]
%so that %the inverse Fourier transform is $\|\widehat f\|_2 = \sqrt{2\pi}\|f\|_2$, $f\in L^2(\RR)$. 
We shall denote by $\mathcal FL^1 = \mathcal FL^1(\RR)$ the Fourier algebra $\mathcal F(L^1(\RR))$. The inverse Fourier transform of a function $h\in\mathcal FL^1(\RR)$ will be denoted by $\check h$ or $\F^{-1}(h)$.

\bd\label{defBS}
Let $(\X,\T)$ be a non-degenerate Banach $L^1(\RR)$-module, and $N$ be a subset of $\X$. The \emph{Beurling spectrum} $\Lambda(N) = \Lambda(N,\T)$ is defined by
\[
\Lambda(N,\T) = \{\l\in\RR: \mbox{ if } f\in L^1 \mbox{ and } fx = 0 \mbox{ for all } x\in N \mbox{ then } \widehat f(\l) = 0 \}.
\] 
\ed

To simplify the notation we shall write $\Lambda(x)$ instead of $\Lambda(\{x\})$, $x\in\X$. 
We refer to \cite[Lemma 3.3]{BK05} for the basic properties of the Beurling spectrum.
In this paper we always assume that there is a $\g > 0$ such that
$\Lambda(\X,\T) \subseteq [-\g, \g].$ As a consequence, our representation $\T$ is uniformly continuous and we have
\[
\T(t) = e^{it\A}, \ t\in\RR,
\] 
for a bounded operator $\A\in B(\X)$ which is the generator of the Banach module $(\X,\T)$ in the usual sense of \cite{BK14, BKU20}.

It is not hard to show that the operators $\T(f)$, $f\in L^1(\RR)$, provide a functional calculus for the generator $\A$. Via the isomorphism of $L^1(\RR)$ and $\mathcal FL^1(\RR)$, we also get the functional calculus $\chT(\widehat f) = \T(f)$, $\widehat f = \mathcal F(f)\in\mathcal FL^1$. In \cite{BK14, BKU20}, this functional calculus was extended to the space $\lloc(\RR)=\{h$: $\RR \to \CC$ such that  
$h\widehat\varphi\in\lhat(\RR)$  for any $\varphi\in L^1(\RR)$
 with 
$\supp{\widehat\varphi}$ compact$\}$. Observe that $\lhat(\RR)\subset\lloc(\RR)$. Moreover,
$\lloc$ is also an algebra under pointwise multiplication.

In the current setting of $\Lambda(\X,\T) \subseteq [-\g, \g]$ extending the functional calculus to $\lloc$ is especially easy.
For $h\in\lloc(\RR)$ we define an operator $\chT(h) \in B(\HH)$ by letting
\begin{equation}\label{cl}
\chT(h)x = h\diamond x : = (h\widehat\varphi)^{\vee} x = \T((h\widehat\varphi)^{\vee})x, \ x\in\X,
\end{equation}
where $\varphi\in L^1(\RR)$ is such that $\supp \widehat\varphi$ is compact and  $\widehat\varphi\equiv 1$ in a neighborhood of $[-\g, \g]$. The vector $\chT(h)x$ is well defined in this way because it is independent of the choice of $\varphi$. 

In this paper, we mostly care about functions $h\in\lloc(\RR)$ 
that satisfy
%from the following definition.
%\bd We say that $h\in\lloc(\RR)$ is an \emph{almost periodic function with a summable Fourier series} if
\begin{equation}\label{aph}
h(\xi) = \sum_{n\in\ZZ} c_n e^{i\xi t_n}, \ \sum_{n\in\ZZ} |c_n| < \infty,\ \xi, t_n\in\RR, n\in\ZZ.
\end{equation}
%The set of all such functions is denoted by $\AP_1$ or $\AP_1(\RR)$.  
%\ed
The key step in the proof of our main results relies on  \cite[Proposition 2.11]{BK14}, which states that for a function $h$ given by \eqref{aph} we have
\begin{equation}
\label{mas}
\chT(h) = \sum_{n\in\ZZ} c_n\T(t_n) \in B(\X).   
\end{equation}

To provide a short proof of Theorem \ref{Bineq}, we also need the following result, which is essentially a restatement of \cite[Theorem 2.14 and Remark 2.8]{BK14}.

\bt\label{prodest}
Let $x\in(\X,\T)$  be such that $\L(x)\subseteq [a,b]$ for some $a < b \in \RR$. Assume also that a function $h\in \lloc(\RR)$ is (real-valued) non-negative, and monotonic %, and convex %\footnote{A continuous function $f:[a,b]\to \RR$ is \emph{convex} if $f(\frac{s+t}2)\le \frac12(f(s)+f(t))$ for all $s,t\in[a,b]$.}
 on $[a,b]$. Additionally, assume that $h$ is convex on $[a,b]$ or satisfies $h^{\prime} \in L^2([a,b])$. 
Then 
\begin{equation}\label{dianorm}
\|h\diamond x\| \le \max\{h(a), h(b)\}\|x\|. 
\end{equation}
\et

We are now ready to provide a proof of Theorem \ref{Bineq}.

\bpf[Proof of Theorem \ref{Bineq}]
According to \eqref{mas}, we have
\[
\T(t)x - x = h\diamond x,
\]
where $h(t) = 1 - e^{i\g t}$. Let $h_1(t) = 1 - \cos{\g t}$ and $h_2(t) = \sin {\g t}$ so that $h = h_1 - ih_2$. Using Theorem \ref{prodest}, we then have
\[
\|\T(t)x - x\| \le \|h_1\diamond x\|+\|h_2\diamond x\| \le 
(1 - \cos{\g t} + \sin{\g t})\|x\|
\]
and the theorem is proved after applying a few standard trigonometric identities.
\epf

To obtain our main results, one cannot simply apply Theorem \ref{prodest}. Nevertheless, the philosophy of the proof remains the same and the approach is rooted in the proofs in \cite{B79, BK14} and, most importantly, in the original proof of M.~Kadec.  

\section{Proofs of the main results.}\label{mainproof}
In this section, $\HH$ is a complex Hilbert space with the $L^1(\RR)$-module structure associated with a representation $\T$, and $\L(\HH,\T)\subseteq [-\pi,\pi]$. It suffices to prove our results in this setting because if the representation $\T$ is such that $\L(\HH,\T)\subseteq [-\g,\g]$, it can be replaced with the representation $\T_\g$ given by $\T_\g(t) = \T(\frac{\pi}{\g} t)$ for which $\L(\HH,\T_\g)\subseteq [-\pi,\pi]$.

For completeness of exposition we recall the standard definitions of Riesz bases and frames.
\bd\label{defRB}
A sequence of vectors $\Phi = \{\varphi_n\}_{n\in\ZZ}$ in $\HH$ forms a \emph{Riesz basis} for $\HH$ with bounds $B>A>0$ if for any sequence $c\in \ell^2(\ZZ)$ one has
\begin{equation}\label{rbineq}
    A\|c\|_2^2 \le \left\|\sum_{n\in\ZZ} c_n\varphi_n \right\|^2
    \le B\|c\|_2^2.
\end{equation}
\ed

\bd\label{defFr}
A sequence of vectors $\Phi = \{\varphi_n\}_{n\in\ZZ}$ in $\HH$ forms a \emph{frame} for $\HH$ with bounds $B>A>0$ if for any vector $x\in \HH$ one has
\begin{equation}\label{frineq}
    A\|x\|^2 \le \sum_{n\in\ZZ} \left|\la x, \varphi_n\ra \right|^2
    \le B\|x\|^2.
\end{equation}
\ed

It is well known that any Riesz basis is a frame with the same bounds. For a frame $\Phi$, we shall denote by $T_\Phi$ its synthesis operator, i.e.~the bounded linear operator from $\ell^2$ to $\HH$ given by $T_\Phi c = \sum_{n\in\ZZ} c_n\varphi_n$, $c = \{c_n\}\subset\ell^2$. Its adjoint operator, $T^*_\Phi: \HH\to\ell^2$, $T^*_\Phi x = \{\la x, \varphi_n\ra\}$ is called the analysis operator. From Definition \ref{defFr} and the standard properties of adjoint operators we see that the operator norms of the analysis and synthesis operators satisfy
\begin{equation}\label{opeq}
    \|T_\Phi\| = \|T^*_\Phi\| \le \sqrt B.
\end{equation}

The starting point for the proof of  Theorem \ref{thmRB}  is the following extension of the Paley-Wiener perturbation lemma.

\begin{lem}[\cite{C95}]\label{PWl}
    Let $T_\Phi$ be the synthesis operator of a frame $\Phi$ with bounds $B\ge A > 0$ and $T_{\widetilde\Phi}$
    be the synthesis operator for a sequence $\widetilde\Phi$. Assume that $\|(T_\Phi - T_{\widetilde\Phi})c\| \le \l 
    \|T_\Phi c\| + \mu\|c\|_2$ for some $\l, \mu \ge 0$ such that $\l+\frac{\mu}{\sqrt{A}} < 1$ and any finitely supported sequence $c \in \ell^2$. Then $\widetilde\Phi$ is a frame
    with bounds $A\left(1-\l- \frac{\mu}{\sqrt{A}}\right)^2$ and $B\left(1+\l+ \frac{\mu}{\sqrt{B}}\right)^2$. If additionally  $\Phi$ is a Riesz basis then so is $\widetilde{\Phi}$.
\end{lem}

We now use the proofs of M.~Kadec \cite{K64} and R.~Balan \cite{Bal97} as a blueprint for obtaining Theorem \ref{thmRB}.

\bpf[Proof of Theorem \ref{thmRB}]
We shall establish applicability of Lemma \ref{PWl} with $\l = 0$ and $\mu = \sqrt{ B}(1- \cos(\pi\delta) +\sin(\pi\delta))$. Observe that in this case we have $\l+\frac{\mu}{\sqrt{A}} < 1$ due to \eqref{deq}.

Let $c \in \ell^2(\mathbb{Z})$ be a finitely supported sequence, and set
    \[
    \begin{split}
       U :&= \left\| (T_\Phi - T_{\widetilde\Phi})c\right\| \\&
= \left\| \sum_{n\in\ZZ}c_n(I - \T(\delta_n))\T(\g_n)x\right\| = \left\| \sum_{n\in\ZZ}c_n h_n\diamond\T(\g_n)x\right\|, 
    \end{split}
\]
where $h_n(t) = 1 - e^{i\delta_n t}$, $\d_n=\widetilde\g_n-\g_n$. Following Kadec, we decompose $h_n$ on $[-\pi, \pi]$  as 
\[
h_n = h_n^1+ h_n^2+h_n^3,
\]
where
\begin{align*}
h_n^1 &\equiv 1 - \frac{\sin{(\pi\delta_n)}}{\pi\delta_n},\\
h_n^2 &= \sum_{\nu = 1}^\infty \frac{(-1)^\nu 2\delta_n \sin(\pi\delta_n)}{\pi(\nu^2-\delta_n^2)}\cos(\nu\cdot),\\ 
h_n^3 &=
i\sum_{\nu = 1}^\infty \frac{(-1)^\nu 2\delta_n \cos(\pi\delta_n)}{\pi((\nu-\frac12)^2-\delta_n^2)}\sin\left(\left(\nu-\frac12\right)\cdot\right).
\end{align*}
Thus, $h_n\diamond \T(\g_n)x = h_n^1\diamond \T(\g_n)x+ h_n^2\diamond \T(\g_n)x+h_n^3\diamond \T(\g_n)x$ and
\[
U \le \left\| \sum_{n\in\ZZ}c_n h_n^1\diamond\T(\g_n)x\right\|+
\left\| \sum_{n\in\ZZ}c_n h_n^2\diamond\T(\g_n)x\right\|+
\left\| \sum_{n\in\ZZ}c_n h_n^3\diamond\T(\g_n)x\right\|.
\]
We proceed to estimate each of the three summands separately. 

For the first one, we use the inequality \eqref{opeq} to get
\[
\begin{split}
    \left\| \sum_{n\in\ZZ}c_n h_n^1\diamond\T(\g_n)x\right\|
&\le \sqrt{B}\left| \sum_{n\in\ZZ}|c_n|^2 \left(1 - \frac{\sin{(\pi\delta_n)}}{\pi\delta_n}\right)^2\right|^{\frac12}, \\&\le \sqrt{B}
\left(1 - \frac{\sin{(\pi\delta)}}{\pi\delta}\right)\|c\|_2
\end{split}
\]

%\[
%\le \sqrt{\frac BA}\left(1 - \frac{\sin{(\pi\delta)}}{\pi\delta}\right)\left\|\sum_{n\in\ZZ}c_n \T(\g_n)x\right\|.
%\]

To estimate the other two terms, we first use
  \cite[Proposition 2.11]{BK14}, see also \eqref{aph} and \eqref{mas}, to conclude that 
\[
h^2_n \diamond \T(\gamma_n)x = \sum_{\nu \in \ZZ\setminus \{0\}} \frac{(-1)^\nu \delta_n \sin(\pi\delta_n)}{\pi(\nu^2-\delta_n^2)}\T({\nu} + \g_n)x 
\]
and 
\[
h_n^3\diamond \T(\gamma_n)x = \sum_{\nu =1}^\infty \frac{(-1)^{\nu} \delta_n \cos(\pi\delta_n)}{\pi((\nu-\frac12)^2-\delta_n^2)}(\T(\nu -\frac{1}{2} + \g_n)+\T(-\nu +\frac{1}{2} + \g_n))x.
\]

Secondly, relying on the facts that the sequence $c$ is finitely supported and the sums over $\nu$ converge absolutely, we  interchange the order of summation and use \eqref{opeq}
%of $\nu$ and $n$ 
to get 
\[
\begin{split}
    \left\| \sum_{n\in\ZZ}c_n h_n^2\diamond\T(\g_n)x\right\| &\le
\left\|\sum_{\nu \in \ZZ\setminus \{0\}} \T(\nu)\left(\sum_{n\in\ZZ} c_n\frac{(-1)^\nu \delta_n \sin(\pi\delta_n)}{\pi(\nu^2-\delta_n^2)}\T(\g_n)x\right)\right\|
\\&
\le \sum_{\nu \in \ZZ\setminus \{0\}} \|\T(\nu)\|\left\|\sum_{n\in\ZZ} c_n\frac{(-1)^\nu \delta_n \sin(\pi\delta_n)}{\pi(\nu^2-\delta_n^2)}\T(\g_n)x\right\|
\\&
\le \sqrt{B} \sum_{\nu \in \ZZ\setminus\{0\}}\left| \sum_{n\in\ZZ}|c_n|^2 \left(\frac{(-1)^\nu \delta_n \sin(\pi\delta_n)}{\pi(\nu^2-\delta_n^2)}\right)^2\right|^{\frac12} 
\\&
\le \sqrt{B} \sum_{\nu \in \ZZ\setminus\{0\}}\frac{\delta \sin(\pi\delta)}{\pi(\nu^2-\delta^2)}\|c\|_2;
\end{split}
\]
%\[\le 
% \sqrt{\frac BA} \sum_{\nu \in \ZZ\setminus\{0\}} \frac{\delta \sin(\pi\delta)}{\pi(\nu^2-\delta^2)}
%\left\|\sum_{n\in\ZZ}c_n \T(\g_n)x\right\|.
%\]
similarly, we obtain
\[
\left\| \sum_{n\in\ZZ}c_n h_n^3\diamond\T(\g_n)x\right\|\le
\sqrt{ B} \sum_{\nu =1}^{\infty} \frac{2 \delta \cos(\pi\delta)}{\pi((\nu-\frac12)^2-\delta^2)}
\left\|c\right\|_2.
\]
%\[\le
% \sqrt{\frac BA} \sum_{\nu =1}^{\infty} \frac{2 \delta \cos(\pi\delta)}{\pi((\nu-\frac12)^2-\delta^2)}
%\left\|\sum_{n\in\ZZ}c_n \T(\g_n)x\right\|.
%\]

Finally, just as in the original proof of Kadec, it  follows that
\begin{equation}\label{fKad}
%    \begin{split}
        U  \le \sqrt{ B}(1- \cos(\pi\delta) +\sin(\pi\delta)) \|c\|_2. 
%        \\& \le \sqrt{\frac BA}(1- \cos(\pi\delta) +\sin(\pi\delta))\left\|\sum_{n\in\ZZ}c_n \T(n)x\right\| %= \l\left\|\sum_{n\in\ZZ}c_n \T(n)x\right\|
%    \end{split}
\end{equation}
%\[
%U \le \sqrt{\frac BA}(1- \cos(\pi\delta) +\sin(\pi\delta))\left\|\sum_{n\in\ZZ}c_n \T(n)x\right\| = 
%\l\left\|\sum_{n\in\ZZ}c_n \T(n)x\right\|
%\]
 Thus, Lemma \ref{PWl} applies and the proof is complete. % since $\l =\sqrt{\frac BA}(1- \cos(\pi\delta) +\sin(\pi\delta)) <1$ due to \eqref{deq}.
\epf

\section{Extension to atomic decompositions}\label{adem1}

The proof in the previous section can be easily adapted for a large class of atomic decompositions in Banach spaces. We follow \cite{CH97} in our exposition of the background in this section. 

\bd\label{adem}
Let  $\X$ be a complex Banach space and $\X'$ be its Banach dual. We say that a pair of sequences $(Y,X)$, $Y =\{y_n\}_{n\in\ZZ}\subset \X'$ and $\X = \{x_n\}_{n\in\ZZ}\subset \X$, forms an \emph{atomic decomposition} of $\X$ with respect to $\ell^p(\ZZ)$ with bounds $A$ and $B$ if the following  conditions hold:
\begin{enumerate}
    \item $\{\la x, y_n\ra\}\in \ell^p$ for each $x\in\X$;
    \item $A\|x\| \le \|\{\la x, y_n\ra\}\|_p \le B\|x\|$ for each $x\in\X$;
    \item $x = \sum_{n\in\ZZ} \la x, y_n\ra x_n$ for each $x\in\X$.
\end{enumerate}
\ed
We refer to linear operators $T_Y: \X\to \ell^p$, $T_Yx = \{\la x, y_n\ra\}$ and $T_X: D(T_X) \subseteq \ell^p
\to \X$, $T_X c = \sum_{n\in\ZZ} c_nx_n$, as the analysis and the synthesis operators of the atomic decomposition $(Y,X)$, respectively. The operator $T_Y$ is automatically bounded by Definition \ref{adem}. The operator $T_X$ need not be bonded in general. However, in this paper we will only consider atomic decompositions for which $T_X$ is a bounded linear operator from $\ell^p$ to $\X$.  %Note that $T_Y\X$ contains all finitely supported sequences for any atomic decomposition.

In place of Lemma \ref{PWl} we then use the following special case of a result by O.~Christensen and C.~Heil (see \cite[Theorem 2]{CH97}).

\bl\label{chlem}
Let $(Y,X)$ be an atomic decomposition for $\X$ with respect to $\ell^p(\ZZ)$, for some $p\in [1,\infty)$, with bounds $A$ and $B$.  Let $W = \{w_n\}
\subset \X$ be such that there exists $\mu\in [0, \frac1B)$ for which $\|(T_X-T_W)c\| \le \mu\|c\|_p$ for any finitely supported sequence $c\in\ell^p$. Then there exists a sequence $Z\subset \X'$ such that $(Z,W)$ is an atomic decomposition for $\X$ with respect to $\ell^p(\ZZ)$ with bounds $A(1+\mu B)^{-1}$ and $B(1-\mu B)^{-1}$. Additionally, $W$ is a basis if and only if $X$ is a basis.
\el

A straightforward adjustment of the arguments in the proofs of the previous section yields the following result.

\bt
Let $\T:\RR\to B(\X)$ be an isometric representation such that $\Lambda(\X,\T) \subseteq [-\g,\g]$ for some $\g > 0$. For  a vector $x\in\X$ and a set $\Gamma = \{\g_n: n\in\ZZ\} \subset\RR$, let $X$ be the sequence $\{\T(\g_n)x\}$. Assume now that a sequence $Y=\{y_n\}_{n\in\ZZ}\subset \X'$ is such that the pair of sequences $(Y, X)$ forms an atomic decomposition for $\X$ with respect to $\ell^p(\ZZ)$, for some $p\in [1,\infty)$, with bounds $B\ge A>0$. Additionally assume that the synthesis operator $T_X: \ell^p\to \X$ is bounded. Let  $\widetilde\Gamma = \{\widetilde\g_n: n\in\ZZ\} \subset\RR$ be such that 
\[
\delta := \sup\{|\widetilde\g_n-\g_n|: n\in\ZZ\} <
 \frac{\pi}{4\g} - \frac1\g\arcsin{\left(\frac{1}{\sqrt2}\left(1-{\frac 1{B\|T_X\|}}\right)\right)}.
\]
Then there exists a sequence $Z\subset \X'$ such that the pair of sequences $(Z, W)$, $W=\{\T(\widetilde\g_n)x\}$, forms an atomic decomposition for $\X$ with respect to $\ell^p(\ZZ)$ with bounds
$A(1+ B\|T_X\|(1-\cos{\delta\g}+\sin{\delta\g}))^{-1}$  and   $B(1-B\|T_X\|(1-\cos{\delta\g}+\sin{\delta\g}))^{-1}$. Moreover, if $X$ is a basis then so is $W$.
\et

\section{Concluding remarks}

One can come up with a multitude of examples where the theorems of this paper can be applied outside of the original realm of exponential bases. It is worth pointing out, however, that our results are already useful for studying frames and Riesz bases of exponential functions defined on sets that are not intervals. In our future research, we expect to address some of the relevant questions pertaining to functions defined on unions of intervals \cite{CL22, PRW24}.

Our interest in the problems discussed in this paper was originally spurred by the questions about dynamical frames posed in \cite{CHFS24} with regard to various dynamical sampling problems \cite{ACMT17, ADK13, AP23}.
Addressing those questions will require considering orbits of group representations that are not isometric. We shall pursue such extensions in our future research as well.

Finally, we note that it would be interesting to investigate if other generalizations of the Kadec $\frac{1}{4}$ theorem, such as the Katsnelson $\frac{1}{4}$ theorem \cite{Ka71} or Avdonin's theorem \cite{Av74}, extend to our setting. This question, however, appears to be extremely challenging due to various complex techniques appearing in their proofs.

\medskip
\noindent {\bf{Acknowledgements}.} Both authors of the paper were supported in part by the NSF grant DMS-2208031. We are thankful to Dr.~Andrei Caragea for helpful discussions.

\bibliographystyle{siam}
\bibliography{refs}

\begin{thebibliography}{10}

\bibitem{ACMT17}
{\sc A.~Aldroubi, C.~Cabrelli, U.~Molter, and S.~Tang}, {\em Dynamical sampling}, Applied and Computational Harmonic Analysis, 42 (2017), pp.~378--401.
\newblock doi: 10.1016/j.acha.2015.08.014.

\bibitem{ADK13}
{\sc A.~Aldroubi, J.~Davis, and I.~Krishtal}, {\em Dynamical sampling: time-space trade-off}, Appl. Comput. Harmon. Anal., 34 (2013), pp.~495--503.

\bibitem{AP23}
{\sc J.~Ashbrock and A.~M. Powell}, {\em Dynamical dual frames with an application to quantization}, Linear Algebra Appl., 658 (2023), pp.~151--185.

\bibitem{Av74}
{\sc S.~A. Avdonin}, {\em On the question of {R}iesz bases of exponential functions in {$L^{2}$}}, Vestnik Leningrad. Univ. Mat. Meh. Astronom.,  (1974), pp.~5--12, 154.

\bibitem{Bal97}
{\sc R.~Balan}, {\em Stability theorems for {F}ourier frames and wavelet {R}iesz bases}, J. Fourier Anal. Appl., 3 (1997), pp.~499--504.

\bibitem{B79}
{\sc A.~G. Baskakov}, {\em Bern\v ste\u\i n-type inequalities in abstract harmonic analysis}, Sibirsk. Mat. Zh., 20 (1979), pp.~942--952, 1164.
\newblock English translation: Siberian Math. J. 20 (1979), no. 5, pp. 665--672 (1980).

\bibitem{BK05}
{\sc A.~G. Baskakov and I.~A. Krishtal}, {\em Harmonic analysis of causal operators and their spectral properties}, Izv. Ross. Akad. Nauk Ser. Mat., 69 (2005), pp.~3--54.
\newblock English translation: Izv. Math. 69 (2005), no. 3, pp. 439--486.

\bibitem{BK14}
\leavevmode\vrule height 2pt depth -1.6pt width 23pt, {\em Memory estimation of inverse operators}, J. Funct. Anal., 267 (2014), pp.~2551--2605.

\bibitem{BKU20}
{\sc A.~G. Baskakov, I.~A. Krishtal, and N.~B. Uskova}, {\em Closed operator functional calculus in {B}anach modules and applications}, J. Math. Anal. Appl., 492 (2020), pp.~124473, 14.

\bibitem{CL22}
{\sc A.~Caragea and D.~G. Lee}, {\em A note on exponential {R}iesz bases}, Sampl. Theory Signal Process. Data Anal., 20 (2022), pp.~Paper No. 13, 14.

\bibitem{C95}
{\sc O.~Christensen}, {\em A {P}aley-{W}iener theorem for frames}, Proc. Amer. Math. Soc., 123 (1995), pp.~2199--2201.

\bibitem{CHFS24}
{\sc O.~Christensen, M.~Hasannasab, F.~M. Philipp, and D.~Stoeva}, {\em The mystery of {C}arleson frames}, Appl. Comput. Harmon. Anal., 72 (2024), pp.~Paper No. 101659, 5.

\bibitem{CH97}
{\sc O.~Christensen and C.~Heil}, {\em Perturbations of {B}anach frames and atomic decompositions}, Math. Nachr., 185 (1997), pp.~33--47.

\bibitem{K64}
{\sc M.~{\u{I}}. Kadec}, {\em The exact value of the {P}aley-{W}iener constant}, Dokl. Akad. Nauk SSSR, 155 (1964), pp.~1253--1254.

\bibitem{Ka71}
{\sc V.~E. Katsnelson}, {\em Bases of exponential functions in {$L^{2}$}}, Funkcional. Anal. i Prilo\v{z}en., 5 (1971), pp.~37--47.

\bibitem{PRW24}
{\sc G.~Pfander, S.~Revay, and D.~Walnut}, {\em Exponential bases for partitions of intervals}, Appl. Comput. Harmon. Anal., 68 (2024), pp.~Paper No. 101607, 22.

\end{thebibliography}
\end{document}